\documentclass{amsart}
\usepackage[applemac]{inputenc}

\usepackage[colorlinks=true,linkcolor=red,citecolor=blue,urlcolor=blue,hypertexnames=false]{hyperref}
\usepackage{bookmark}
\usepackage{amsthm,thmtools,amssymb,amsmath,amscd}

\usepackage{fancyhdr}
\usepackage{esint}
\usepackage{enumerate}

\usepackage{pictexwd,dcpic}

\swapnumbers
\declaretheorem[name=Theorem,numberwithin=section]{thm}
\declaretheorem[name=Remark,style=remark,sibling=thm]{rem}
\declaretheorem[name=Lemma,sibling=thm]{lemma}

\declaretheorem[name=Corollary,sibling=thm]{cor}

\numberwithin{equation}{section}

\usepackage{cleveref}
\crefname{lemma}{Lemma}{Lemmata}
\crefname{prop}{Proposition}{Propositions}
\crefname{thm}{Theorem}{Theorems}
\crefname{cor}{Corollary}{Corollaries}
\crefname{defn}{Definition}{Definitions}
\crefname{example}{Example}{Examples}
\crefname{rem}{Remark}{Remarks}
\crefname{assum}{Assumption}{Assumptions}
\crefname{nota}{Notation}{Notation}
\usepackage{autonum}

\newcommand{\ti}{\tilde}
\newcommand{\wt}{\widetilde}

\newcommand{\cn}{\colon}
\newcommand{\sub}{\subset}

\newcommand{\mr}{\mathring}

\newcommand{\bbR}{\mathbb{R}}

\newcommand{\bbS}{\mathbb{S}}
\newcommand{\bbH}{\mathbb{H}}

\newcommand{\bbM}{\mathbb{M}}

\newcommand{\8}{\infty}

\newcommand{\al}{\alpha}
\newcommand{\be}{\beta}
\newcommand{\g}{\gamma}
\newcommand{\de}{\delta}
\newcommand{\e}{\epsilon}
\newcommand{\ka}{\kappa}
\newcommand{\la}{\lambda}

\newcommand{\om}{\omega}
\newcommand{\s}{\sigma}
\newcommand{\Si}{\Sigma}
\newcommand{\p}{\varphi}

\newcommand{\vt}{\vartheta}

\newcommand{\D}{\Delta}
\newcommand{\G}{\Gamma}

\newcommand{\cL}{\mathcal{L}}

\newcommand{\cS}{\mathcal{S}}

\newcommand{\del}{\partial}
\newcommand{\n}{\nabla}
\newcommand{\II}[2]{\mrm{II}\br{#1,#2}}

\newcommand{\rt}{\sqrt}

\newcommand{\ip}[2]{\left\langle #1,#2 \right\rangle}
\newcommand{\fr}[2]{\frac{#1}{#2}}
\newcommand{\x}{\times}

\DeclareMathOperator{\graph}{graph}

\DeclareMathOperator{\vol}{vol}

\newcommand{\pf}[1]{\begin{proof}{\parskip\baselineskip{ #1}} \end{proof}}
\newcommand{\eq}[1]{\begin{equation}\begin{alignedat}{2} #1 \end{alignedat}\end{equation}}

\newcommand{\br}[1]{\left(#1\right)}
\newcommand{\abs}[1]{\lvert #1\rvert}
\newcommand{\enum}[1]{\begin{enumerate}[(i)] #1 \end{enumerate}}

\newcommand{\ra}{\rightarrow}

\newcommand{\mt}{\mapsto}

\newcommand{\mrm}{\mathrm}

\newcommand{\hp}{\hphantom}
\newcommand{\q}{\quad}

\begin{document}

\title{The Minkowski inequality in de Sitter space}
\date{\today}
\thanks{Funded by the "Deutsche Forschungsgemeinschaft" (DFG, German research foundation); Project "Quermassintegral preserving local curvature flows"; Grant number SCHE 1879/3-1.}
\author[J. Scheuer]{Julian Scheuer}
\address{Department of Mathematics, Columbia University
New York, NY 10027, USA}
\email{\href{mailto:jss2291@columbia.edu}{jss2291@columbia.edu}; \href{mailto:julian.scheuer@math.uni-freiburg.de}{julian.scheuer@math.uni-freiburg.de}}

\urladdr{\href{https://home.mathematik.uni-freiburg.de/scheuer/}{https://home.mathematik.uni-freiburg.de/scheuer/}}

\begin{abstract}		
The classical Minkowski inequality in the Euclidean space provides a lower bound on the total mean curvature of a hypersurface in terms of the surface area, which is optimal on round spheres. In this paper we employ a locally constrained inverse mean curvature flow to prove a properly defined analogue in the Lorentzian de Sitter space. 
\end{abstract}
\maketitle

\section{Introduction}

\subsection*{The main results}
In this paper we prove an optimal Minkowski inequality for compact, spacelike and mean convex hypersurfaces in the upper branch of the $(n+1)$-dimensional Lorentzian de Sitter space. To state the main result, we provide the involved terminology. Let $n\geq 2$ and denote by $\bbM^{n+2}_{1}$ the $(n+2)$-dimensional Minkowski space, i.e. $\bbR^{n+2}$ with metric
\eq{\ip{v}{w}=-v^{0}w^{0}+\sum_{\al=1}^{n+1}v^{\al}w^{\al}.}
We define the upper branch of the de Sitter space by
\eq{\bar\bbS^{n+1}_{1}=\{y\in \bbM^{n+2}_{1}\cn \ip{y}{y}=1,\q y^{0}>0\}}
and understand this submanifold to be equipped with the induced metric. Then $\bar{\bbS}^{n+1}_{1}$ is a Lorentzian manifold with constant sectional curvature 1, see \Cref{DeSitter} for a detailed discussion. 
The following theorem is the main result of this paper.

\begin{thm}\label{main}
Let $\Si\sub \bar{\bbS}^{n+1}_{1}$ be a spacelike, compact, connected and mean-convex hypersurface. Then there holds
\eq{\int_{\Si}H_{1}\leq \vol(\hat\Si)+\p(\abs{\Si})}
with equality precisely if $\Si$ is totally umbilic. Here $\p\cn \bbR_{+}\ra \bbR$ is the strictly increasing function which gives equality on the $y^{0}$-slices.
\end{thm}

We review the standard terminology used here very briefly and refer to \Cref{DeSitter} for a more detailed discussion. A hypersurface $\Si\sub \bar{\bbS}^{n+1}_{1}$ is called spacelike, if its induced metric is Riemannian. $\Si$ is called mean-convex, if for a suitable unit normal vector field $\nu$ the normalized mean curvature $H_{1}$ with respect to $-\nu$ is positive.
$\hat \Si$ denotes the region enclosed by the slice $\{y^{0}=0\}$ and $\Si$ and $\vol(\hat\Si)$ its enclosed volume, see \eqref{vol}. Finally, $\abs{\Si}$ denotes the surface area of $\Si$. 

In order to prove \Cref{main}, we employ a locally constrained inverse mean curvature flow. This flow is designed to preserve the surface area and to increase the quantity
\eq{W_{2}(\Si)=\int_{\Si}H_{1}-\vol(\hat\Si).}
We will prove that it converges smoothly to a coordinate slice, which will imply \Cref{main}. To state the result about the curvature flow, we introduce some more notation. The space $\bar{\bbS}^{n+1}_{1}$ is isometric to the space $\bbR_{+}\x \bbS^{n}$ with the warped product metric
\eq{\bar g=-dr^{2}+\vt^{2}(r)\s,}
where $\s$ is the round metric on $\bbS^{n}$ and $\vt=\cosh,$ cf. \Cref{cosh}. For a hypersurface $\Si$ we define
\eq{u=-\bar g(\vt\del_{r},\nu)}
to be the support function, where $\nu$ is the future directed normal vector field on $\Si$, where the time orientation is inherited from $\bbM^{n+2}_{1}$.

\begin{thm}\label{flow-main}
Let $\Si\sub \bar{\bbS}^{n+1}_{1}$ be a spacelike, compact, connected and mean-convex hypersurface. Then there exists a unique immortal solution
\eq{x\cn [0,\8)\times \bbS^{n}&\ra \bar{\bbS}_{1}^{n+1}}
which satisfies
\eq{\label{flow-main-1}\dot{x}&=\br{u-\fr{\vt'}{H_{1}}}\nu\\
		x(0,\bbS^{n})&=x_{0},}
where $x_{0}$ is an embedding of $\Si$.
As $t\ra \8$, the embeddings $x(t,\cdot)$ converge smoothly to a coordinate slice, which is uniquely determined by $\abs{\Si}$.
\end{thm}

\subsection*{Background}

\subsubsection*{Minkowski inequality}
For a closed and convex hypersurface in the Euclidean space $\bbR^{n+1}$, the Minkowski inequality states that
\eq{\label{Minkowski}\int_{\Si}H_{1}\geq \abs{\bbS^{n}}^{\fr 1n}\abs{\Si}^{\fr{n-1}{n}}}
with equality precisely when $\Si$ is a round sphere.
For surfaces $n=2$ this was originally proved by Minkowski \cite{Minkowski:/1903}. Using inverse curvature flows, \eqref{Minkowski} was, among more general estimates called {\it{Alexandrov-Fenchel inequalities}}, generalized to starshaped and mean-convex hypersurfaces in \cite{GuanLi:08/2009}. Using Huisken's and Ilmanen's weak inverse mean curvature flow \cite{HuiskenIlmanen:/2001} one can replace the starshapedness by outward minimality. It is open until today whether \eqref{Minkowski} holds for general mean-convex hypersurfaces. The Michael-Simon-Sobolev inequality \cite{MichaelSimon:05/1973} gives an estimate on the total mean curvature, however its optimal constant is not the desired optimal constant in \eqref{Minkowski}, also compare \cite{Brendle:07/2019}. In case $n=2$, there is an $L^{2}$-stability result for immersed surfaces \cite{KuwertScheuer:06/2019}, namely there holds
			\eq{\left|\fr{1}{\rt{\abs{\Si}}}\int_{\Si}H_{1} - 2\rt{\pi} \right|  \leq  c\|\mr A \|_{L^{2}(\Si)}^{2},}
where $\mr A$ is the trace-free part of the second fundamental form. Also see \cite{DalphinHenrotMasnouTakahashi:10/2016} for a comprehensive overview over related results for the case of closed hypersurfaces of the Euclidean space. For hypersurfaces with free boundary on a cone there are similar results \cite{Cruz:03/2019}.
Both sides of \eqref{Minkowski} can be considered as special cases of the quermassintegrals $W_{k}(\Si)$. We refer to \cite{GallegoSolanes:/2005,GaoHugSchneider:/2001,Solanes:/2006} for a comprehensive introduction and useful further references. In the Euclidean case \eqref{Minkowski} then reads
\eq{W_{2}^{\bbR^{n+1}}(\Si)\geq c_{n}W_{1}^{\bbR^{n+1}}(\Si)^{\fr{n-1}{n}}.}
The superscript is added to indicate the special structure of $W_{k}(\Si)$ in the Euclidean case.

There are analogues of these quantities in the hyperbolic and spherical space $\bbH^{n+1}$ resp. $\bbS^{n+1}$. However, the $W_{k}(\Si)$ are then not given by curvature integrals, but by linear combinations of them. While the quantities $W_{1}(\Si)$ are up to dimensional constants given by surface area in all of the three spaceforms, the values of $W_{2}(\Si)$ in the non-flat spaces are (up to dimensional constants)
\eq{W_{2}^{\bbH^{n+1}}(\Si)=\int_{\Si}H_{1}-\vol(\hat \Si),\q W_{2}^{\bbS^{n+1}}(\Si)=\int_{\Si}H_{1}+\vol(\hat\Si),}
where $\hat \Si$ is the region enclosed by $\Si$.

A possible generalization of \eqref{Minkowski} to other ambient manifolds would then be to derive an inequality between $W_{2}(\Si)$ and $W_{1}(\Si)$ (we drop the superscript again for the following informal discussion).

In the hyperbolic space, \cite{GallegoSolanes:/2005} treats the case of convex hypersurfaces through a rough estimate which does not characterize the case of equality, while in \cite{WangXia:07/2014} an optimal inequality was proved in the class of horospherically convex hypersurfaces. Variants of such estimates involving total mean curvature, surface area and possibly other quantities, are contained for example in \cite{BorisenkoMiquel:/1999,De-LimaGirao:04/2016,Natario:08/2015,WeiXiong:/2015}.

In the sphere there are also many variants of \eqref{Minkowski}, for example \cite{GiraoPinheiro:12/2017,MakowskiScheuer:11/2016,WeiXiong:/2015} and in other ambient spaces there are weighted Minkowski-type inequalities \cite{BrendleHungWang:01/2016,GeWangWuXia:03/2015,McCormick:09/2018,ScheuerXia:11/2019,Wang:/2015,Wei:04/2018,Xia:/2016a}.

In the de Sitter space we are only aware of one similar result, which recently appeared \cite{AndrewsHuLi:03/2019}. The Alexandrov-Fenchel inequality proved in this paper is deduced as a corollary from its dual version in the hyperbolic space with the help of a well-known duality method available for convex hypersurfaces of the sphere and hyperbolic/de Sitter space, see \cite{GaoHugSchneider:/2001,Gerhardt:/2015} and \Cref{DeSitter}.

Note that the inequality in \Cref{main} can not be deduced by duality from hyperbolic space, as it is supposed to be shown for non-convex hypersurfaces, for which the duality is not defined.

\subsubsection*{Locally constrained curvature flows}

The strategy to prove \Cref{main} is to employ a suitably defined curvature flow, i.e. a variation of the hypersurface $\Si$ which is defined through its curvature and possibly lower order quantities. This method has become very popular for the deduction of geometric inequalities during the past decades. For example, in \cite{Huisken:/1987} Huisken studied the {\it{volume preserving mean curvature flow}}, where the mean curvature flow is modified by addition of a global term, namely the averaged mean curvature. Long-time existence and convergence to a round sphere is proven, if the initial datum is convex. This reproved the isoperimetric inequality in all dimensions for convex domains. Similar nonlocal flows have been widely used to prove other geometric inequalities in the Euclidean and hyperbolic space \cite{Andrews:/2001,Athanassenas:/1997,Cabezas-RivasMiquel:/2007,IvakiStancu:/2013,McCoy:/2003,McCoy:02/2004,McCoy:/2005,Sinestrari:10/2015}. 
These nonlocal flows are hard to study due to the nonlocal term involved and the results mentioned above usually required preservation of convexity at the least.

A new kind of volume preserving curvature flow was invented by Pengfei Guan and Junfang Li in \cite{GuanLi:/2015}. In the Euclidean space it is
\eq{\label{GL}\dot{x}=\fr{1}{2n}\D_{\Si}\abs{x}^{2}=(1-uH_{1})\nu,}
where $H_{1}$ is the normalized mean curvature and $u$ the support function. This flow obviously preserves the enclosed volume and it can be calculated that it decreases the surface area unless $\Si$ is a round sphere. This gives a further proof of the isoperimetric inequality. The major advantage over the classical volume preserving mean curvature flow is that it preserves the starshapedness and hence the technical obstructions are lower. Hence there has been some recent interest in these so-called {\it{locally constrained curvature flows}}. Later, the flow \eqref{GL} was also studied in other ambient manifolds \cite{GuanLiWang:/2019}. From that idea, it was tempting to modify \eqref{GL} in such a way that the modified flow preserves other geometric quantities. For example, the flow
\eq{\label{GL2}\dot{x}=\br{\fr{1}{H_{1}}-u}\nu}
in $\bbR^{n+1}$ preserves the surface area and decreases the total mean curvature, giving another proof of the Minkowski inequality for starshaped and mean-convex hypersurfaces. To quickly complete the list of previous literature on these locally constrained flows, we refer to the related results in \cite{GuanLi:/2018,ScheuerWangXia:11/2018,ScheuerXia:11/2019}

Note that while \eqref{GL2} is very easy to treat in  $\bbR^{n+1}$ as it is basically a rescaling of the inverse mean curvature flow originally treated in \cite{Gerhardt:/1990,Urbas:/1990}, it seems hard to study the proper modification of this flow in the other spaceforms. In the hyperbolic space there is a partial unpublished result \cite{BrendleGuanLi:/} which requires an additional initial gradient smallness assumption, while in the sphere there is no result. 

The purpose of this paper is to prove the full convergence result \Cref{flow-main} of the correct version of \eqref{GL2} in the Lorentzian de Sitter space under the most natural assumptions, which are spacelikeness and mean convexity, and in turn obtain the Minkowski inequality for such hypersurfaces.

\subsection*{Outline}
In \Cref{DeSitter} we spend some time to review the geometry of de Sitter space and to introduce our notation. We take some care here, as we will also deduce some maybe not so well known relations between hypersurfaces of de Sitter space and their duals in the hyperbolic space. In particular we will deduce the dual flow of \eqref{flow-main-1} in the hyperbolic space. In \Cref{Ev} and \Cref{Apriori} we collect the relevant evolution equations and deduce the required a priori estimates to obtain long-time existence. In \Cref{proof} we prove that the flow converges to a coordinate slice and complete the proof of the Minkowski inequality.

\section{Geometry of de Sitter space and duality}\label{DeSitter}

We recall some facts about hypersurfaces of semi-Riemannian manifolds as well as basic properties of the de Sitter space, most of which can be found in \cite{ONeill:/1983}.

For a semi-Riemannian manifold $(\bar M,\ip{\cdot}{\cdot})$ with Levi-Civita connection $\bar D$ and a hypersurface $M$ we have the Gaussian formula for vector fields $V,W$ on $M$ (which are smoothly extended to $\bar M$),
\eq{\bar D_{V}W=D_{V}W+\II{V}{W},}
where $D$ is the Levi-Civita connection of the metric induced by the inclusion $\iota\cn M\ra \bar{M}$ and the decomposition is orthogonal. 
The normal part $\II{\cdot}{\cdot}$ is called the vector valued second fundamental form.
There holds the Gauss equation \cite[p.~100]{ONeill:/1983},
\eq{\label{Gauss}\ip{R(V,W)X}{Y}&=\ip{\bar R(V,W)X}{Y}+\ip{\II{V}{X}}{\II{W}{Y}}\\
			&\hp{=}-\ip{\II{V}{Y}}{\II{W}{X}},}
where we used the curvature tensor convention from \cite{ONeill:/1983}, 
\eq{R(X,Y)Z=D_{Y}D_{X}Z-D_{X}D_{Y}Z-D_{[Y,X]}Z}
for all vector fields $X,Y,Z$ on $M$. 
The shape operator $S$ of $M$ derived from a normal  $N$, which is defined by
\eq{\ip{S(V)}{W}=\ip{\II{V}{W}}{N},}
satisfies the Weingarten equation \cite[p.~107]{ONeill:/1983}
\eq{\label{Weingarten}S(V)=-\bar D_{V}N.}

Let $n\geq 2$ and $\bbM^{n+2}_{1}$ be the $(n+2)$-dimensional Minkowski space with metric
\eq{\ip{v}{w}=-v^{0}w^{0}+\sum_{\al=1}^{n+1}v^{\al}w^{\al}.}
The Lorentzian de Sitter space is defined as the hyperquadric
\eq{\bbS^{n+1}_{1}=\{y\in \bbM^{n+2}_{1}\cn \ip{y}{y}=1\}.}

Differentiating $1=\ip{\g}{\g}$ along an arbitrary curve $\g$ in $\bbS^{n+1}_{1}$, we obtain that the normal space of $\bbS^{n+1}_{1}\sub \bbM^{n+2}_{1}$ is spanned by the spacelike position vector field $y$ and thus $\bbS^{n+1}_{1}$ has sign $1$ in $\bbM^{n+2}_{1}$. From \eqref{Weingarten} we also obtain that the shape operator of $\bbS^{n+1}_{1}$ is
\eq{S(V)=-\bar D_{V}y=-V}
and hence $\II{\cdot}{\cdot}=-\ip{\cdot}{\cdot}y$.
We obtain from \eqref{Gauss} that the Riemann tensor $R$ of $\bbS^{n+1}_{1}$ satisfies
\eq{\ip{R(V,W)X}{Y}=\ip{V}{X}\ip{W}{Y}-\ip{V}{Y}\ip{W}{X}.}
It follows that $\bbS^{n+1}_{1}$ has constant sectional curvature
\eq{K(e_{i},e_{j})=\fr{\ip{R(e_{i},e_{j})e_{i}}{e_{j}}}{\ip{e_{i}}{e_{i}}\ip{e_{j}}{e_{j}}}=\fr{\ip{e_{i}}{e_{i}}\ip{e_{j}}{e_{j}}}{\ip{e_{i}}{e_{i}}\ip{e_{j}}{e_{j}}}=1}
for every orthogonal unit vectors $e_{i},e_{j}$, \cite[p.~77]{ONeill:/1983}. 

For the calculations in this paper it will be convenient to have a warped product structure for $\bbS^{n+1}_{1}$ and hence we prove the following lemma.

\begin{lemma}\label{cosh}
The Lorentzian manifold $\bbS^{n+1}_{1}$ is isometric to the warped product $\bbR\x \bbS^{n}$
with metric
\eq{\bar g=-dr^{2}+\vt^{2}(r)\s,}
where $\s$ is the round metric on $\bbS^{n}$ and $\vt=\cosh$.
The hyperbolic space 
\eq{\bbH^{n+1}:=\{\ti y\in \bbM_{1}^{n+2}\cn \ip{\ti y}{\ti y}=-1,\q \ti y^{0}>0\}} is diffeomorphic to $\bbR_{+}\x \bbS^{n}$ with metric
\eq{\ti{\bar g}=d\ti r^{2}+\ti \vt ^{2}(\ti r)\s,}
where $\ti \vt=\sinh$. 
\end{lemma}

\pf{
It follows from \cite[p.~111]{ONeill:/1983}, that the map
\eq{\phi\cn \bbR\x \bbS^{n}&\ra \bbS^{n+1}_{1}\sub \bbM^{n+2}_{1}\\
			(r,p)&\mt \br{\vt'(r),\vt(r)p}}
is a diffeomorphism. Let $x^{i}$ be local coordinates on $\bbS^{n}$. The pullback metric of $\phi$ is calculated as follows:
\eq{\bar g(\del_{r},\del_{r})=\left|\br{\vt,\vt' p}\right|^{2}=-1,}
\eq{\bar g(\del_{r},\del_{x^{i}})=0,}
since $p$ is normal to $\del_{x^{i}}p$ in $\bbR^{n+1}$ and
\eq{\bar g(\del_{x^{i}},\del_{x^{j}})=\vt^{2}(r)\s_{ij}}
by definition. 
The hyperbolic case is proven by replacing $\vt$ by $\ti \vt$ in the above proof.
}

In order to employ the duality between $\bar{\bbS}^{n+1}_{1}$ and $\bbH^{n+1}$ we will need the following formulae. The first follows from the previous proof and the second is a straightforward computation of the push-forward of $\vt\del_{r}$. 

\begin{lemma}\label{ambient}
As quantities of the ambient manifold, the function $\vt'(r)$ is given by
\eq{\vt'(r)=-\ip{y}{e}}
and the vector field $\vt\del_{r}$ is 
\eq{\vt\del_{r}=e-\ip{e}{y}y,}
where we identified $(r,p)$ with $y\in \bbS^{n+1}_{1}\sub \bbR^{n+2}_{1}$ and $e=e_{0}$.

Similarly, the function $\ti \vt(\ti r)$ is given by
\eq{\ti\vt'(\ti{r})=-\ip{\ti y}{e}}
while the vector field $\ti\vt\del_{\ti r}$ is given by
\eq{\ti\vt\del_{\ti r}=-e-\ip{e}{\ti y}\ti y.}
\end{lemma}

\begin{rem}
We will see later, that the curvature flow we will employ to prove the main theorem will only be parabolic if $\vt'>0$. Hence will will restrict to the upper branch $\{r>0\}\sub\bbS^{n+1}_{1}$.
\end{rem}

 \subsection*{Spacelike hypersurfaces of de Sitter space}
 
Let $\Si\sub\bbS^{n+1}_{1}$ be a spacelike, compact, connected hypersurface. The manifold $\bbS^{n+1}_{1}$ is globally hyperbolic \cite[Thm.~1.4.2]{Gerhardt:/2006} and 
\eq{\cS_{0}=\{0\}\x \bbS^{n}}
a Cauchy hypersurface, terminologies which are defined e.g. in \cite[Def.~1.3.7, 1.3.8]{Gerhardt:/2006}. From \cite[Prop.~1.6.3, Rem.~1.6.4]{Gerhardt:/2006} we obtain that $\Si$ is a smooth graph over $\cS_{0}$, i.e. in the coordinates from \Cref{cosh} we have
\eq{\Si=\{(\rho(x^{i}),x^{i})\cn (x^{i})\in \cS_{0} \},}
where in the sequel latin indices range between $1$ and $n$, while greek indices range from $0$ to $n$. As mentioned above, we will assume $\rho>0$. In order to shorten notation we will often write 
\eq{x^{0}=r}
and hence a point $x\in \bbS^{n+1}_{1}$ will have the coordinate representation $x=(x^{\al})$ in a given coordinate system.

Now we describe the geometry of $\Si$ in terms of the ambient geometry and $\rho$, starting with the introduction of some notation. Fix a local coordinate system $(\xi^{i})$ for $\Si$. For tensors on $\Si$ we will use the coordinate based notation, e.g. the induced metric $g$ will be written as $g=(g_{ij})$, where
\eq{g_{ij}=g(\del_{\xi^{i}},\del_{\xi^{j}}).}
If $\n$ denotes the Levi-Civita connection of $g$, in order to shorten the appearance of evolution equations, we will denote the coordinate functions of covariant derivatives of tensors by the use of semi-colons:
If $A=(a^{i_{1}\dots i_{k}}_{j_{1}\dots j_{l}})$ is a $k$-contravariant and $l$-covariant tensor (or merely a function), its covariant derivative is denoted by
\eq{\n_{\del_{x^{m}}} A=\br{a^{i_{1}\dots i_{k}}_{j_{1}\dots j_{l};m}}.}

For spacelike hypersurfaces $\Si$ we define $\nu$ to be the future directed (timelike) normal, i.e.
\eq{\bar g(\del_{r},\nu)<0.}
We derive the shape operator $S=(h^{i}_{j})$ of $\Si$ from $-\nu$,
\eq{h_{ij}:=g_{ik}h^{k}_{j}=-\bar g\br{\II{\del_{\xi^{i}}}{\del_{\xi^{j}}},\nu}}
and call the tensor $(h_{ij})$ the second fundamental form of $\Si$. The eigenvalues $(\ka_{i})_{1\leq i\leq n}$ of the shape operator are called the principal curvatures of $\Si$. This definition implies that
\eq{\II{\del_{\xi^{i}}}{\del_{\xi^{j}}}=h_{ij}\nu.}
This is in accordance with the convention in \cite[Thm.~1.1.2]{Gerhardt:/2006} and we use this reference for further formulae.

From \cite[(1.6.13)]{Gerhardt:/2006} we obtain that the second fundamental form of the slices $\{x^{0}=r\}$ is given by
\eq{\bar{h}_{ij}:=\fr{\vt'(r)}{\vt(r)}\bar g_{ij}.}
The induced metric in terms of $\rho$ is
\eq{g_{ij}=-\rho_{;i}\rho_{;j}+\vt^{2}(\rho)\s_{ij}=-\rho_{;i}\rho_{;j}+\bar g_{ij}.}
Defining
\eq{v^{2}=1-\vt^{-2}\s^{ij}\rho_{;i}\rho_{;j},}
the second fundamental form satisfies
\eq{\label{h}v^{-1}h_{ij}=\rho_{;ij}+\bar h_{ij},}
\cite[(1.6.11)]{Gerhardt:/2006}, where we note that in this reference the past directed normal is used. Also note that in order to control the property of a hypersurface of being spacelike, one has to ensure $v^{2}>0$.

The function 
\eq{u:=\fr{\vt}{v}=-\bar g(\vt\del_{r},{\nu})}
is of special interest and can be regarded as a generalized support function.

\subsection*{Graphical hypersurfaces of hyperbolic space}
Let $\ti\Si\sub \bbH^{n+1}$ be starshaped with respect to the origin in the coordinates given by \Cref{cosh},
\eq{\ti\Si=\{(\ti\rho(x^{i}),x^{i})\cn (x^{i})\in \cS_{0}\}.}
With the corresponding notation as for the de Sitter space, we always chose the normal $\ti\nu$ to $\Si$ to point in the same direction as $\del_{\ti r}$, i.e.
\eq{\ti{\bar g}(\del_{\ti r},\ti\nu)>0.}
We derive the shape operator of $\Si$ from $-\ti\nu$ and obtain from \cite[(1.5.10)]{Gerhardt:/2006} that
\eq{\label{h-hyp}\ti h_{ij}\ti v^{-1}=-\ti\rho_{;ij}+\ti{\bar h}_{ij},}
where
\eq{\ti v^{2}=1+\ti \vt^{-2}\s^{ij}\ti\rho_{;i}\ti\rho_{;j}.}
We also may define the support function of $\ti \Si$ to be
\eq{\ti u=\ti{\bar g}(\ti\vt\del_{\ti r},\ti \nu).}

\subsection*{Duality}
There is an important relation between strictly convex hypersurfaces of the hyperbolic space containing the origin and strictly convex hypersurfaces of the upper branch of de Sitter space. It states that these two sets are in one-to-one relation through the idempotent Gauss maps
\eq{\ti x:=\nu,\q x:=\ti \nu.}
Furthermore the respective principal curvatures satisfy the relation
\eq{\ti \ka_{i}=\ka_{i}^{-1}}
and the induced metric of the dual hypersurfaces are given by
\eq{\label{dual-h}\ti g_{ij}=h_{ik}h^{k}_{j},\q g_{ij}=\ti h_{ik}\ti h^{k}_{j}.}
We refer to \cite[Thm.~10.4.4, 10.4.5, 10.4.9]{Gerhardt:/2006} for a more thorough discussion. In the theory of curvature flows, the method of duality has successfully be employed several times, e.g. \cite{BIS4,Gerhardt:/2015,Yu:/2017}.

Later we will employ the following relations between the support functions and the respective height functions of the dual hypersurfaces.

\begin{lemma}\label{dual}
Let $\Si\sub\bar{\bbS}^{n+1}_{1}$ be a strictly convex, compact and spacelike hypersurface. Then its dual $\ti\Si\sub\bbH^{n+1}$ is starshaped with respect to the origin and the height/support functions satisfy
\eq{\vt'(r)=\ti u,\q \ti\vt'(\ti r)=u.} 
\end{lemma}

\pf{We use \Cref{ambient} and $\ip{x}{\ti x}=0$. There holds
\eq{u=-\ip{e-\ip{e}{x}x}{\ti x}=-\ip{e}{\ti x}=\ti\vt'(\ti r)}
and 
\eq{\ti u=-\ip{e+\ip{e}{\ti x}\ti x}{x}=-\ip{e}{x}=\vt'(r).}
}

\subsection*{Volume, surface area and Minkowski identities}
A function $f$ on a spacelike hypersurface $\Si\sub\bbS^{n+1}_{1}$ is in $L^{1}(\Si)$, if the differential $n$-form $\abs{f}d\om_{g}$ has a finite integral over $\Si$, in which case we write
\eq{\int_{\Si}f:=\int_{\bbS^{n}}f\,d\om_{g}.}
Here $d\om_{g}$ is the Riemannian volume form on $\Si$. 
For a spacelike hypersurface $\Si=\graph \rho$ as above we define the {\it{enclosed volume}} as in \cite[Sec.~4]{Makowski:01/2013} by
\eq{\label{vol}\vol(\hat\Si)=\int_{\bbS^{n}}\int_{0}^{\rho(\cdot)}\fr{\rt{\det(\bar g_{ij}(s,\cdot))}}{\rt{\det(\s_{ij})}}\,ds}
and the surface area by
\eq{\abs{\Si}=\int_{\Si}1.}

Now we recall some crucial crucial integral identities, which are known as Min\-kowski identities in the Riemannian context. We restrict to those we need, while more general versions were proved in \cite{Kwong:01/2016}.
First of all let
\eq{H_{k}(\ka)=\fr{1}{\binom{n}{k}}\sum_{1\leq i_{1}\leq\dots\leq i_{k}\leq n}\ka_{i_{1}}\dots\ka_{i_{k}},\q 1\leq k\leq n,}
be the normalized elementary symmetric polynomials, where we also define $H_{0}=1$.
We calculate 
\eq{\label{rho}\vt'(\rho)_{;ij}=\vt' \rho_{;i}\rho_{;j}+\vt \rho_{;ij}=\vt'\rho_{;i}\rho_{;j}+\fr{\vt}{v}h_{ij}-\vt'\bar g_{ij}=-\vt' g_{ij}+\fr{\vt}{v}h_{ij}.}
Regarding the functions $H_{k}=H_{k}(g_{ij},h_{ij})$ in dependence of the metric and the second fundamental form, in space forms the tensors
\eq{H_{k}^{ij}=\fr{\del H_{k}}{\del h_{ij}}}
are divergence free. Taking into account that
\eq{g_{ij}H_{k}^{ij}=kH_{k-1}, \q 1\leq k\leq n,}
we obtain from tracing \eqref{rho} with respect to $H_{k}^{ij}$ and integration
\eq{\label{Mink}\int_{\Si}\vt'H_{k-1}=\int_{\Si}u H_{k},\q 1\leq k\leq n.}
In the hyperbolic space the same relations can be deduced from \eqref{h-hyp}.

\section{Evolution equations}\label{Ev}

To prove the Minkowski inequality for $\Si$, we will make use of normal variations and therefore we recall some known variational formulae and deduce the ones we specifically need here.

A quite general treatment also appears in \cite[Ch.~2]{Gerhardt:/2006}, while our framework does not fit into the setting discussed there. Hence we take some more care in this section.

\subsection*{General evolution equations}
Let $T>0$, $M$ compact and connected and
\eq{x\cn (0,T)\x M\ra \bar{\bbS}^{n+1}_{1}}
a time-dependent family of embeddings of the spacelike hypersurfaces 
\eq{\Si_{t}=x(t,M).}
Let us denote their normal velocity by $f$, i.e.
\eq{\label{var}\dot{x}=f\nu.}

From \cite[p.~94]{Gerhardt:/2006} we get
\eq{\del_{t}g_{ij}=2fh_{ij},}
and
\eq{\dot{\nu}=f_{;k}g^{kj}x_{;j}}
where $(g^{kj})$ is the inverse of $g$ and a dot denotes the covariant time derivative of a tensor field along the curve $x(\cdot,\xi)$.

As \cite{Gerhardt:/2006} uses a different convention for the Riemann tensor, for the reader's convenience we deduce the evolution of the shape operator here.

\begin{lemma}\label{gen-ev-h}
There holds
\eq{\del_{t}h^{j}_{i}={f_{;i}}^{j}-fh^{j}_{k}h^{k}_{i}+f\de^{j}_{i}.}
\end{lemma}

\pf{
Let $\bar \n$ denote the Levi-Civita connection of $\bar g$. We differentiate the Weingarten equation
\eq{\bar \n_{x_{;i}}\nu=h^{k}_{i}x_{;k}}
covariantly with respect to time and get
\eq{\del_{t}h^{k}_{i}x_{;k}+h^{k}_{i}\dot{x}_{;k}&=\bar\n_{\dot x}\bar\n_{x_{;i}}\nu=\bar\n_{x_{;i}}\bar\n_{\dot x}\nu+\bar{R}(x_{;i},\dot x)\nu=\bar\n_{x_{;i}}\bar\n_{\dot x}\nu+fx_{;i}}
and hence, after multiplying with $x_{;j}$,
\eq{g_{kj}\del_{t}h^{k}_{i}&=fg_{ij}-fh^{k}_{i}h_{kj}+\bar g(\bar\n_{x_{;i}}(f_{;k}g^{kl}x_{;l}),x_{;j})\\
					&=fg_{ij}-fh^{k}_{i}h_{kj}+f_{;ij}.}
}

\begin{lemma}\label{GeomQ}
Volume, surface area and total mean curvature evolve by
\eq{\del_{t}\vol(\hat\Si_{t})=\int_{\Si_{t}}f,\q\del_{t}\abs{\Si_{t}}=n\int_{\Si_{t}}fH_{1}}
and
\eq{\del_{t}\int_{\Si_{t}}H_{1}=(n-1)\int_{\Si_{t}}fH_{2}+\del_{t}\vol(\hat\Si_{t}).}
\end{lemma}

\pf{
According to \cite[(2.4.21)]{Gerhardt:/2006}, the radial function $\rho$ satisfies
\eq{\fr{\del\rho}{\del t}=fv,}
where we note again the flip of our normal compared to that reference.
Hence
\eq{\del_{t}\vol(\hat\Si_{t})=\int_{\bbS^{n}}fv\fr{\rt{\det(\bar g_{ij}(\rho(\cdot),\cdot))}}{\rt{\det(\s_{ij})}}=\int_{\bbS^{n}}f\fr{\rt{\det(g_{ij})}}{\rt{\det(\s_{ij})}}=\int_{\Si}f.}

Thee holds
\eq{\del_{t}\rt{\det (g_{ij})}=\fr{\det(g_{ij})g^{kl}\del_{t}g_{kl}}{2\rt{\det (g_{ij})}}=nfH_{1}\rt{\det(g_{ij})}.}
The second claim follows.
We calculate
\eq{\del_{t}\int_{\Si_{t}}H_{1}&=\fr 1n\int_{\Si_{t}}f(n^{2}H_{1}^{2}-\abs{S}^{2}+n)\\
				&=\fr{2}{n}\binom{n}{2}\int_{\Si_{t}}fH_{2}+\int_{\Si_{t}}f.}
}

\subsection*{Specific evolution equations}

In order to prove the Minkowski inequality, we want to build a flow that preserves the surface area and increases the quantity
\eq{W_{2}(\Si)=\int_{\Si}H_{1}-\vol(\hat\Si).}
A natural flow to consider is the following locally constrained inverse mean curvature flow:
\eq{\label{Flow}\dot{x}=\br{u-\fr{\vt'(\rho)}{H_{1}}}\nu,}
which indeed has the desired properties:

\begin{lemma}\label{Pres}
Along \eqref{Flow} the surface area is preserved and the quantity
\eq{W_{2}(\Si_{t})=\int_{\Si_{t}}H_{1}-\vol(\hat\Si_{t})}
is non-decreasing and strictly increasing unless the flow hypersurfaces are totally umbilic. 
\end{lemma}

\pf{
From \Cref{GeomQ} and \eqref{Mink} we obtain
\eq{\fr 1n \del_{t}\abs{\Si_{t}}=\int_{\Si_{t}}(uH_{1}-\vt')=0}
and
\eq{\del_{t}W_{2}(\Si_{t})=(n-1)\int_{\Si_{t}}\br{uH_{2}-\fr{\vt'H_{2}}{H_{1}}}=(n-1)\int_{\Si_{t}}\br{\vt'H_{1}-\fr{\vt'H_{2}}{H_{1}}}\geq 0,}
where we used
\eq{H_{2}=H_{1}^{2}-\|\mr A\|^{2},}
where $\mr A$ is the trace-free part of the second fundamental form.
}

We need the evolution equation for the radial and support function and therefore define
\eq{\cL=\del_{t}-\fr{\vt'}{nH_{1}^{2}}\D-\vt{\rho_{;}}^{k}\del_{k}.}

\begin{lemma}
\enum{
\item The radial function $\rho$ satisfies
\eq{\label{Ev-rho}\cL\rho=\vt-\fr{2\vt'}{H_{1}}v^{-1}+\fr{\vt'^{2}}{\vt H_{1}^{2}}+\fr{\vt'^{2}}{n\vt H_{1}^{2}}\|\n \rho\|^{2}.}
\item The support function $u$ satisfies
\eq{\label{Ev-u}\cL u=-\fr{\vt'}{nH_{1}^{2}}\|\mr A\|^{2}u-\fr{\vt^{2}}{H_{1}}\|\n \rho\|^{2}.}

}

\end{lemma}

\pf{
(i)~The choice of the normal implies
\eq{\del_{t}\rho=\br{u-\fr{\vt'}{H_{1}}}v^{-1}.}
From \eqref{h} we obtain
\eq{\del_{t}\rho-\fr{\vt'}{nH_{1}^{2}}\D\rho&=\br{u-\fr{\vt'}{H_{1}}}v^{-1}-\fr{\vt'}{nH_{1}^{2}}g^{ij}(v^{-1}h_{ij}-\bar h_{ij})\\
			&=\br{u-\fr{\vt'}{H_{1}}}v^{-1}-\fr{\vt'}{H_{1}}v^{-1}+\fr{\vt'^{2}}{n\vt H_{1}^{2}}g^{ij}\bar g_{ij}\\
			&=\vt v^{-2}-\fr{2\vt'}{H_{1}}v^{-1}+\fr{\vt'^{2}}{\vt H_{1}^{2}}+\fr{\vt'^{2}}{n\vt H_{1}^{2}}\|\n \rho\|^{2}.}
The result follows from
\eq{\|\n\rho\|^{2}=\br{\bar g^{ij}+\fr{\bar g^{ik}\rho_{;k}\bar g^{jl}\rho_{;l}}{v^{2}}}\rho_{;i}\rho_{;j}=1-v^{2}+\fr{(1-v^{2})^{2}}{v^{2}}=\fr{1-v^{2}}{v^{2}}.}
			
(ii)~We calculate, that the vector field $\vt\del_{x^{0}}$ is conformal:
\eq{(\vt\del_{x^{0}})_{;\al}=\vt'r_{;\al}\del_{x^{0}}+\vt\bar\G^{\be}_{0\al}\del_{x^{\be}}.}
The Christoffel-symbols of $\bar g$ are
\eq{\bar \G^{\be}_{0\al}&=\fr 12 \bar g^{\be\de}\br{\fr{\del}{\del x^{0}}\bar g_{\al\de}+\fr{\del}{\del x^{\al}}\bar g_{0\de}-\fr{\del}{\del x^{\de}}\bar g_{0\al}}\\
			\\
			&=\begin{cases} 0,\q &\be=0\\
						 \fr{\vt'}{\vt}\de^{i}_{\al},\q &\be=i
				\end{cases}
			}
and hence
\eq{(\vt\del_{x^{0}})_{;\al}=\vt'r_{;\al}\del_{x^{0}}+\vt'\de^{i}_{\al}\del_{x^{i}}=\vt'\de^{0}_{\al}\del_{x^{0}}+\vt'\de^{i}_{\al}\del_{x^{i}}=\vt'\de^{\be}_{\al}\del_{x^{\be}}=\vt'\del_{x^{\al}}.}
Hence
\eq{\del_{t}u=-\bar g(\bar \n_{\dot x}(\vt\del_{r}),\nu)-\bar g(\vt\del_{r},\dot\nu)=f\vt'-\bar g(\vt\del_{r},x_{;j})f_{;k}g^{kj},}
\eq{u_{;i}=-h^{k}_{i}\bar g(\vt\del_{r},x_{;k})}
and
\eq{\label{u}u_{;ij}=-h^{k}_{i;j}\bar g(\vt\del_{r},x_{;k})-\vt' h_{ij}+h^{k}_{i}h_{kj}u.}
Hence
\eq{\del_{t}u-\fr{\vt'}{nH_{1}^{2}}\D u&=-\fr{\vt'}{nH_{1}^{2}}\|A\|^{2}u+\vt'\br{u-\fr{\vt'}{H_{1}}}+\fr{\vt'^{2}}{H_{1}}\\
				&\hp{=}-\bar g(\vt\del_{r},x_{;j})\br{u-\fr{\vt'}{H_{1}}}_{;k}g^{kj}+\fr{\vt'}{H_{1}^{2}}{H_{1;}}^{k}\bar g(\vt\del_{r},x_{;k})\\
				&=-\fr{\vt'}{nH_{1}^{2}}\br{\|A\|^{2}-nH_{1}^{2}}u-\bar g(\vt\del_{r},x_{;j})u_{;k}g^{kj}+\bar g(\vt\del_{r},x_{;j})\fr{\vt}{H_{1}}{\rho_{;}}^{j}\\
				&=-\fr{\vt'}{nH_{1}^{2}}\|\mr A\|^{2}u+\vt g(\n \rho,\n u)-\fr{\vt^{2}}{H_{1}}\|\n \rho\|^{2} .}
}

We conclude this section with the full evolution of the shape operator.

\begin{lemma}
There holds
\eq{\label{Ev-h}\cL h^{i}_{j}&=-\vt'\br{\fr{\|A\|^{2}}{nH_{1}^{2}}h^{j}_{i}-\fr{2}{H_{1}}h^{m}_{i}h^{j}_{m}+h^{j}_{i}}+u\br{\de^{j}_{i}-\fr{h^{j}_{i}}{H_{1}}}+\fr{\vt'}{H_{1}^{2}}(H_{1}\de^{j}_{i}-h^{j}_{i})\\
				&\hp{=}+\fr{\vt}{H_{1}^{2}}{H_{1;}}^{j}\rho_{;i}+\fr{\vt}{H_{1}^{2}}H_{1;i}{\rho_{;}}^{j}-2\fr{\vt'}{H_{1}^{3}}H_{1;i}{H_{1;}}^{j}.}\end{lemma}

\pf{
From \Cref{gen-ev-h} we have
\eq{\del_{t}h^{j}_{i}={\br{u-\fr{\vt'}{H_{1}}}_{;i}}^{j}-\br{u-\fr{\vt'}{H_{1}}}h^{j}_{k}h^{k}_{i}+\br{u-\fr{\vt'}{H_{1}}}\de^{j}_{i}.}
We have to expand the second order term. There holds
\eq{\br{u-\fr{\vt'}{H_{1}}}_{;i}=u_{;i}-\fr{\vt}{H_{1}}\rho_{;i}+\fr{\vt'}{H_{1}^{2}}H_{1;i}}
and
\eq{\br{u-\fr{\vt'}{H_{1}}}_{;ij}&=u_{;ij}-\fr{\vt'}{H_{1}}\rho_{;i}\rho_{;j}+\fr{\vt}{H_{1}^{2}}H_{1;j}\rho_{;i}-\fr{\vt}{H_{1}}\rho_{;ij}+\fr{\vt}{H_{1}^{2}}H_{1;i}\rho_{;j}\\
				&\hp{=}-2\fr{\vt'}{H_{1}^{3}}H_{1;i}H_{1;j}+\fr{\vt'}{H_{1}^{2}}H_{1;ij}.}
In order to replace the term $H_{1;ij}$, we have to use the Codazzi and Gauss equation \eqref{Gauss}.
As in \cite[p.~76]{ONeill:/1983} we define
\eq{R(\del_{k},\del_{l})\del_{j}=R^{m}_{jkl}\del_{m}.}
 There holds
\eq{nH_{1;ij}=h^{k}_{k;ij}=h^{k}_{i;kj}=h^{k}_{i;jk}+h^{m}_{i}R^{k}_{mkj}-h^{k}_{m}R^{m}_{ikj}={h_{ij;k}}^{k}+h^{m}_{i}R^{k}_{mkj}-h^{k}_{m}R^{m}_{ikj}.}
We use \eqref{Gauss} and $\bar g(\nu,\nu)=-1$ to deduce

\eq{R^{i}_{jkl}=g(R(\del_{k},\del_{l})\del_{j},\del_{m})g^{im}=\de^{i}_{l}g_{jk}-\de^{i}_{k}g_{jl}-h^{i}_{l}h_{jk}+h^{i}_{k}h_{jl}}
and hence
\eq{nH_{1;ij}&={h_{ij;k}}^{k}+h^{m}_{i}R^{k}_{mkj}-h^{k}_{m}R^{m}_{ikj}\\
			&={h_{ij;k}}^{k}+h^{m}_{i}(\de^{k}_{j}g_{mk}-ng_{mj}-h^{k}_{j}h_{mk}+h^{k}_{k}h_{mj})\\
			&\hp{=}-h^{k}_{m}(\de^{m}_{j}g_{ik}-\de^{m}_{k}g_{ij}-h^{m}_{j}h_{ik}+h^{m}_{k}h_{ij})\\
			&={h_{ij;k}}^{k}-(n-1)h_{ij}-h^{m}_{i}h_{mk}h^{k}_{j}+nH_{1}h^{m}_{i}h_{mj}\\
			&\hp{=}-h_{ij}+nH_{1}g_{ij}+h^{k}_{m}h^{m}_{j}h_{ik}-\|A\|^{2}h_{ij}\\
			&={h_{ij;k}}^{k}-nh_{ij}+nH_{1}h^{m}_{i}h_{mj}+nH_{1}g_{ij}-\|A\|^{2}h_{ij}.}
We obtain from \eqref{h} and \eqref{h}:
\eq{\br{u-\fr{\vt'}{H_{1}}}_{;ij}&=-h^{k}_{i;j}\bar g(\vt\del_{r},x_{;k})-\vt'h_{ij}+h^{k}_{i}h_{kj}u-\fr{\vt'}{H_{1}}\rho_{;i}\rho_{;j}+\fr{\vt}{H_{1}^{2}}H_{1;j}\rho_{;i}\\
				&\hp{=}-\fr{\vt}{H_{1}}(v^{-1}h_{ij}-\bar h_{ij})+\fr{\vt}{H_{1}^{2}}H_{1;i}\rho_{;j}-2\fr{\vt'}{H_{1}^{3}}H_{1;i}H_{1;j}\\
				&\hp{=}+\fr{\vt'}{nH_{1}^{2}}({h_{ij;k}}^{k}-nh_{ij}+nH_{1}h^{m}_{i}h_{mj}+nH_{1}g_{ij}-\|A\|^{2}h_{ij})\\
				&=-h^{k}_{i;j}\bar g(\vt\del_{r},x_{;k})-\vt'h_{ij}+h^{k}_{i}h_{kj}u+\fr{\vt}{H_{1}^{2}}H_{1;j}\rho_{;i}\\
				&\hp{=}-\fr{u}{H_{1}}h_{ij}+\fr{2\vt'}{H_{1}} g_{ij}+\fr{\vt}{H_{1}^{2}}H_{1;i}\rho_{;j}-2\fr{\vt'}{H_{1}^{3}}H_{1;i}H_{1;j}\\
				&\hp{=}+\fr{\vt'}{nH_{1}^{2}}{h_{ij;k}}^{k}-\fr{\vt'}{H_{1}^{2}}h_{ij}+\fr{\vt'}{H_{1}}h^{m}_{i}h_{mj}-\fr{\vt'}{nH_{1}^{2}}\|A\|^{2}h_{ij}.}
Thus
\eq{\cL h^{i}_{j}&=-\vt'h_{i}^{j}+\fr{\vt}{H_{1}^{2}}{H_{1;}}^{j}\rho_{;i}-\fr{u}{H_{1}}h_{i}^{j}+\fr{\vt'}{H_{1}} \de^{j}_{i}+\fr{\vt}{H_{1}^{2}}H_{1;i}{\rho_{;}}^{j}-2\fr{\vt'}{H_{1}^{3}}H_{1;i}{H_{1;}}^{j}\\
				&\hp{=}-\fr{\vt'}{H_{1}^{2}}h_{i}^{j}+2\fr{\vt'}{H_{1}}h^{m}_{i}h_{m}^{j}-\fr{\vt'}{nH_{1}^{2}}\|A\|^{2}h_{i}^{j}+u\de^{j}_{i}\\
				&=-\vt'\br{\fr{\|A\|^{2}}{nH_{1}^{2}}h^{j}_{i}-\fr{2}{H_{1}}h^{m}_{i}h^{j}_{m}+h^{j}_{i}}+u\br{\de^{j}_{i}-\fr{h^{j}_{i}}{H_{1}}}+\fr{\vt'}{H_{1}^{2}}(H_{1}\de^{j}_{i}-h^{j}_{i})\\
				&\hp{=}+\fr{\vt}{H_{1}^{2}}{H_{1;}}^{j}\rho_{;i}+\fr{\vt}{H_{1}^{2}}H_{1;i}{\rho_{;}}^{j}-2\fr{\vt'}{H_{1}^{3}}H_{1;i}{H_{1;}}^{j}.}
}

\section{A priori estimates}\label{Apriori}
We establish $C^{2}$-estimates for \eqref{Flow} and use standard regularity theory for parabolic equations to conclude smooth convergence of the flow to a round sphere.

We will assume throughout this section that $\Si\sub \bar{\bbS}^{n+1}_{1}$ is a smooth, closed, connected, spacelike and mean-convex hypersurface.

Then the differential operator $\cL$ is strictly parabolic at $\Si$ and hence the flow \eqref{Flow} has a unique solution for a short time $T^{*}$ with initial hypersurface $\Si_{0}=\Si$. The a priori estimates of this section refer to this solution.

\subsection*{Estimates up to first order}

\begin{lemma}
Along the flow \eqref{Flow} there holds for all $(t,\xi)\in (0,T^{*})\x M$:
\enum{
\item
\eq{\label{barrier}\min_{M}\rho(0,\cdot)\leq \rho(t,\xi)\leq \max_{M}\rho(0,\cdot)}
\item\eq{u(t,\xi)\leq \max_{M}u(0,\cdot).}
}
\end{lemma}

\pf{
(i)~The radial function $\rho$ satisfies
\eq{\del_{t}\rho=\br{u-\fr{\vt'}{H_{1}}}v^{-1}.}
From \eqref{h} we obtain that spatial maximal points of $\rho$ we have
\eq{0\geq \D\rho = nH_{1}-n\fr{\vt'}{\vt}}
and hence
\eq{u-\fr{\vt'}{H_{1}}\leq 0.}
Hence $\max \rho$ is non-increasing and the reverse estimate proves that $\min\rho$ is non-decreasing.

(ii)~Directly follows from \eqref{Ev-u} and the maximum principle.
}

\subsection*{Curvature estimates}

Tracing \eqref{Ev-h} yields the evolution of the normalized mean curvature:
\eq{\cL H_{1}=\fr{\vt'}{nH_{1}}\br{\|A\|^{2}-nH_{1}^{2}}+\fr{2\vt}{nH_{1}^{2}}g(\n H_{1},\n \rho)-\fr{2\vt'}{nH_{1}^{3}}\|\n H_{1}\|^{2}.}

\begin{lemma}
\enum{
\item There holds for all $(t,\xi)\in (0,T^{*})\x M$:
\eq{H_{1}(t,\xi)\geq \min_{M}H_{1}(0,\cdot).}
\item There exists a constant $c=c(M_{0})$ such that
\eq{H_{1}\leq c.}
}
\end{lemma}

\pf{
(i)~Follows directly from the maximum principle with the help of.
(ii)~Define
\eq{w=\log H_{1}+\la u-\rho,}
where $\la$ is determined to be a large number. Then
\eq{\cL w&=\fr{\cL H_{1}}{H_{1}}+\fr{\vt'}{nH_{1}^{2}}\|\n \log H_{1}\|^{2}+\la \cL u-\cL\rho\\
		&=\fr{\vt'}{nH_{1}^{2}}\|\mr A\|^{2}(1-\la u)+\fr{2\vt}{nH_{1}^{2}}g(\n \log H_{1},\n\rho)-\fr{\vt'}{nH_{1}^{2}}\|\n \log H_{1}\|^{2}-\fr{\la \vt^{2}}{H_{1}}\|\n \rho\|^{2}\\
		&\hp{=}-\vt+\fr{2\vt'}{H_{1}}v^{-1}-\fr{\vt'^{2}}{\vt H_{1}^{2}}-\fr{\vt'^{2}}{n\vt H_{1}^{2}}\|\n \rho\|^{2}\\
		&\leq -\vt +\fr{2\vt'}{H_{1}}v^{-1}+\fr{c}{H_{1}^{2}}\|\n \rho\|^{2}}
for large $\la$.
If $H_{1}$ is too large, this is negative and hence we obtain a bound on $H_{1}$ by the maximum principle. 
}

\begin{lemma}
There exists a constant $c=c(M_{0})$ such that
\eq{\|A\|^{2}\leq c.}
\end{lemma}

\pf{We estimate the largest principal curvature $\ka_{n}$. By a well known trick, cf. \cite[p.~500]{Gerhardt:11/2011}, it suffices to estimate the evolution of $h^{n}_{n}=\ka_{n}$. By \eqref{Ev-h} and with the help of all previously deduced bounds we obtain constants $\e$ and $c$ such that
\eq{\cL h^{n}_{n}\leq -\e \ka_{n}^{3} + c\ka_{n}^{2} + c +c\|\n H_{1}\|\|\n\rho\|-2\e\|\n H_{1}\|^{2}<0, }
if $\ka_{n}$ is sufficiently large. This threshold depends on $\e$ and $c$. The proof is complete.
}

\begin{cor}
The flow \eqref{Flow} starting from $\Si$ exists for all times and satisfies uniform estimates in $C^{m}(\bbS^{n})$ for all $m\geq 0$.
\end{cor}

\pf{
After the previously established $C^{2}$-estimates this is standard from parabolic regularity \cite{Krylov:/1987} applied to the graph function $\rho$.
}

\section{Completion of the proof}\label{proof}

In order to complete the proof, we have to show that the flow converges to a round sphere. 

\begin{lemma}
The flow \eqref{Flow} converges to a uniquely determined coordinate slice and hence \Cref{flow-main} holds.
\end{lemma}

\pf{

Along the flow the quantity
\eq{W_{2}(\Si_{t})=\int_{\Si_{t}}H_{1}-\vol(\hat\Si_{t})}
is clearly bounded and non-decreasing.
Hence
\eq{\del_{t}W_{2}(\Si_{t})=(n-1)\int_{\Si_{t}}\fr{\vt'\|\mr A\|^{2}}{H_{1}}\ra 0}
as $t\ra \8$ and thus every subsequential $C^{\8}$-limit $\Si_{\8}$ must be totally umbilical. As $H_{1}>0$, $\Si_{\8}$ is strictly convex. The dual hypersurface $\ti\Si_{\8}$, which is given by the Gauss map
\eq{\ti x=\nu\cn M\ra \bbH^{n+1}\sub\bbM^{n+2}_{1},}
is thus also totally umbilical, cf. \cite[Thm.~10.4.4]{Gerhardt:/2006} and hence a geodesic sphere. From \eqref{dual-h} we obtain
\eq{\abs{\Si_{\8}}=\int_{\ti\Si_{\8}}\ti H_{n},}
which is, up to an additive constant, the $(n-1)$-quermassintegral $\wt W_{n-1}(\ti\Si_{\8})$ in $\bbH^{n+1}$. See \cite{WangXia:07/2014} for a definition. As $\abs{\Si_{\8}}$ is independent of the subsequential limit, so is $\wt{W}_{n-1}(\ti\Si_{\8})$ and hence the radius and the principal curvatures of $\ti\Si_{\8}$ are uniquely determined. This implies that $H_{1}$ takes the same value $0<H_{*}<1$ for every subsequential limiting hypersurface $\Si_{\8}$.
In particular, the function $H_{1}$ converges uniformly to $H_{*}$ along the flow.
Let
\eq{r_{0}=\tanh^{-1}(H_{*}).}
We claim that the flow converges to the slice $\{r=r_{0}\}$ smoothly and prove this in several steps.
\enum{
\item For every subsequential limit $\Si_{\8}$ with corresponding radial function $\rho_{\8}$ there holds 
\eq{\max_{M}\rho_{\8}\geq r_{0},}
 since at a point where the maximum is attained we have
 \eq{0\geq \D \rho_{\8}=nH_{*}-n\tanh(\rho_{\8})}
 and hence $\rho_{\8}\geq r_{0}$ due to the monotonicity of $\tanh$.

\item Define 
\eq{\p(t)=\max_{M}\rho(t,\cdot)-r_{0}=\rho(t,\xi_{t})-r_{0}.}
Then $\p$ is Lipschitz and hence differentiable almost everywhere. Let $\e>0$ and fix some $T_{\e}>0$ to be specified later.
 Let $t>T_{\e}$ be a point of differentiability of $\p$ and suppose 
\eq{\p(t)\geq \e.}
Note that this condition is only non-void for bounded $\e\leq \e_{0}$, due to the barrier estimates.
Then there holds
\eq{\p'(t)=\del_{t}\rho(t,\xi_{t})=\br{u-\fr{\vt'}{H_{1}}}v^{-1}(t,\xi_{t})=\vt-\fr{\vt'}{H_{1}},}
where the right-hand side is evaluated at $\rho(t,\xi_{t})$. 
 We estimate
 \eq{\p'(t)=\fr{\vt'}{H_{1}}\br{\coth(\rho(t,\xi_{t}))H_{1}-1}&\leq \fr{\vt'}{H_{1}}\br{\coth(r_{0}+\e)H_{*}-1}+c_{r_{0}}\abs{H_{1}-H_{*}}\\
 				&=-\fr{\vt'}{H_{1}}\fr{\vt'(r_{0})}{\vt(r_{0})\vt'^{2}(\eta)}\e+c_{r_{0}}\abs{H_{1}-H_{*}}\\
				&=\br{-\fr{\vt'}{H_{1}}\fr{\vt'(r_{0})}{\vt(r_{0})\vt'^{2}(\eta)}+\fr{c_{r_{0}}}{\e}\abs{H_{1}-H_{*}}}\e,}
where we have used the mean value theorem and $\eta\in[r_{0},r_{0}+\e_{0}]$. As $H_{1}\ra H_{*}$, we may now choose $T=T_{\e}$, which only depends on $\e$ and the initial data, such that for all $t>T_{\e}$ with $\p(t)\geq \e$ there holds
\eq{\p'(t)\leq -\de_{\e}.}
From \cite[Lemma~4.2]{Scheuer:05/2015} it follows that
\eq{\limsup_{t\ra \8} \max_{M}\rho(t,\cdot)\leq r_{0}.}
\item Combining (i) and (ii) we obtain
\eq{\lim_{t\ra \8}\max_{M}\rho(t,\cdot)=r_{0}.}
\item A similar argument applied to $\min \rho$ implies 
\eq{\lim_{t\ra \8}\min_{M}\rho(t,\cdot)=r_{0}.} 
}
Hence the unique limit is the slice $\{r=r_{0}\}$ and the proof is complete.
}

\pf{
Now we prove \Cref{main}. So let $\Si\sub \bar{\bbS}^{n+1}_{1}$ be a spacelike, compact, connected and mean-convex hypersurface. According to \Cref{flow-main}, we may deform $\Si$ in an infinite amount of time to a coordinate slice $S_{r_{\8}}$ of radius $r_{\8}$.
Define 
\eq{\p_{1}(r)=\abs{\{x^{0}=r\}},\q \p_{2}(r)=W_{2}(\{x^{0}=r\}).}
From \Cref{GeomQ} we see that both of these function are strictly increasing functions of $r$.
 By \Cref{Pres} there holds
\eq{W_{2}(\Si)\leq W_{2}(S_{r_{\8}})=\p_{2}(r_{\8})=\p_{2}\circ \p_{1}^{-1}(\p_{1}(r_{\8}))=\p_{2}\circ \p_{1}^{-1}(\abs{\Si}).}
From the proof of \Cref{Pres} we obtain that if we have equality in this inequality, $\Si$ must be totally umbilic, otherwise the flow would increase $W_{2}$ strictly.

It remains to show that total umbilicity implies equality. So let $\Si$ be totally umbilic. Now it suffices to prove that this property is preserved along the flow \eqref{Flow}, since we know that this flow will deform $\Si$ into a coordinate slice on which equality holds. Furthermore, if all flow hypersurfaces are totally umbilic, $W_{2}$ is constant. As $\abs{\Si_{t}}$ is constant anyway, the equality must already hold on $\Si$.

So let us prove that \eqref{Flow} preserves total umbilicity. If $\Si$ is totally umbilic, it must be strictly convex. Hence for a short time the flow hypersurfaces $\Si_{t}$ are strictly convex. As in \cite[Sec.~4]{BIS4} and employing \Cref{dual} we can calculate, that up to a tangential diffeomorphism the dual hypersurfaces $\ti\Si_{t}\sub\bbH^{n+1}$ satisfy the flow equation
\eq{\dot{\ti x}=\br{u-\fr{\vt'}{H_{1}(\ka_{i})}}\ti\nu=\br{\ti\vt'-n\ti u\fr{\s_{n}(\ti\ka)}{\s_{n-1}(\ti \ka)}}\ti \nu,}
where $\s_{k}$ is the $k$-th elementary symmetric polynomial. This is a locally constrained curvature flow of contracting type in hyperbolic space. In the Euclidean space such kind of flows were studied in \cite{GuanLi:/2018}. Although in the non flat spaces no satisfactory convergence results are available in general, we can still show that this flow preserves the total umbilicity. It has the property that it preserves the $(n-1)$-quermassintegral $\wt{W}_{n-1}(\ti\Si_{t})$ in $\bbH^{n+1}$, while it decreases $\wt{W}_{n}(\ti\Si_{t})$ and we have
\eq{\label{dual-flow-pres}\del_{t}\wt{W}_{n}(\ti\Si_{t})<0}
at $t>0$ unless $\ti\Si_{t}$ is totally umbilic. Since by assumption $\ti \Si$ is totally umbilic, we have
\eq{\wt{W}_{n}(\ti\Si)=\psi(\wt{W}_{n-1}(\ti\Si)),}
cf. \cite[Thm.~1.1]{WangXia:07/2014} with a suitable function $\psi$. If $\ti\Si_{t}$ was not totally umbilic for some $t>0$, by \eqref{dual-flow-pres} we would obtain 
\eq{\wt{W}_{n}(\ti\Si_{t})<\psi(\wt{W}_{n-1}(\ti\Si_{t})),}
which violates \cite[Thm.~1.1]{WangXia:07/2014}. Hence the $\ti\Si_{t}$ must be totally umbilic and, going back to the flow in de Sitter space, $\Si_{t}$ must be totally umbilic. Hence \eqref{Flow} preserves the total umbilicity and the proof is complete.
}

\section*{Acknowledgments}
This work was made possible through a research scholarship the author received from the DFG and which was carried out at Columbia University in New York. JS would like to thank the DFG, Columbia University and especially Prof.~Simon Brendle for their support.

Parts of this work were written during a visit of the author to McGill University in Montréal. The author would like to thank McGill and Prof.~Pengfei Guan for their hospitality and helpful discussions.

\bibliographystyle{shamsplain}
\bibliography{Bibliography}
\end{document}